\documentclass[english,russian]{amsart}
\usepackage{amssymb,amsmath}

\newtheorem{thm}{Theorem}[section]
\newtheorem{crl}[thm]{Corollary}
\newtheorem{prp}[thm]{Proposition}
\newtheorem{lm}[thm]{Lemma}

\newcommand{\der}{\partial}

\newcommand{\eps}{\varepsilon}

\newcommand{\sudda}[1]{}

\begin{document}
\title{ Weak Leibniz algebras and Transposed Poisson algebras }
\author{A.S. Dzhumadil'daev}

\address
{Institute of Mathematica, Pushkin street 125,  Almaty,  050000,
Kazakhstan} \email{dzhuma@hotmail.com}


\maketitle

\begin{abstract}
An algebra with identities $[a,b]c=2a(bc)-2b(ac), a[b,c]=2(ab)c-2(ac)b$ is called weak Leibniz. We show that weak Leibniz operad is self-dual and is not Koszul. We establish that polarization of  any weak Leibniz algebra is transposed Poisson, and, conversely, polarization of any transposed Poisson algebra is weak Leibniz. 
\end{abstract}

\section{Introduction}
%A main field $K$ is supposed to be field of characteristic $p\ne 2,$ if otherwise is not stated. 

\subsection{Leibniz algberas and weak Leibniz algebras} An algebra $A$ is called {\it weak Leibniz}, if it satisfies the following identities
$$[a,b]c=2\,a(bc)-2\,b(ac),$$
$$a[b,c]=2\,(ab)c-2\,(ac)b,$$
for any $a,b,c\in A.$ Here 
$$[a,b]=ab-ba$$
is Lie commutator.  

{\bf Example.} An algebra is called (two-sided) Leibniz \cite{Loday}, if 
$$(ab)c=a(bc)-b(ac),$$
$$a(bc)=(ab)c-(ac)b,$$
are identities. Any Leibniz algebra is weak Leibniz.

Proof repeats arguments of J.L. Loday.
By left Leibniz identity 
$$(ab)c=a(bc)-b(ac),$$
$$(ba)c=b(ac)-a(bc).$$
Addition of these relations gives us Loday's observation: anti-commutator of any elements of left-Leibniz algebra is left central,
$$\{a,b\}c=0.$$
If we use substraction of these relations instead of addition we obtain left weak Leibniz identity 
$$[a,b]c=2(a(bc)-b(ac)).$$

Similarly, by right Leibniz identity
$$a(bc)=(ab)c-(ac)b,$$
$$a(cb)=(ac)b-(ab)c,$$
and addition of these relations show that anti-commutator of any elements of right-Leibniz algebra is  right-central, 
$$a\{b,c\}=0.$$
If we use substraction instead of addition we obtain right weak Leibniz identity 
$$a[b,c]=2((ab)c-(ac)b)).$$

In particular any Lie algebra is weak Leibniz. Important property of Leibniz agebras. {\it If $L$ is simple Leibniz, then $L$ is Lie. } Proof is easy. 
Since anti-commutators generate an ideal, $\{a,b\}=0$ for any $a,b\in L.$ This means that $L$ is Lie. 

{\bf Example.}  Let $U$ be associative commutative algebra with derivation $\der:U\rightarrow U$ and for any $u,v\in U$ endow $U$ by new multiplication 
$\times =\times_{u,v} $ constructed by  
$$a\times b=u(a\der(b)-b(\der(a))+v a b.$$
Then an algebra with vector space $U$ and multiplication $\times$ is weak Leibniz. In general this algebra is not Leibniz.

{\bf Example.} Let $S$ be finite set of indices and $\eps_s\in {\bf R}$ for any $s\in S.$ An infinite dimensional algebra $A(S)$ with base elements $e_i, i\in {\bf Z}$ and multiplication
$$e_i\times e_j =(j-i)e_{i+j}+\sum_{s\in S} \eps_s e_{i+j+s}$$
is weak Leibniz. It is simple. Proof of this fact is easy also. 

Denote by  $A(S)^-$ an algebra with base $\{e_i, i\in {\bf Z}\}$ and multiplication 
$[a,b]=1/2(a\times b-b\times a).$ Then $[e_i,e_j]=(j-i)e_{i+j},$ and the algebra $A(S)^-$ is isomorphic to Witt algebra (= Lie algebra of formal vector fields on the circle). If $A(S)$ has an ideal $J$ then $J$ should be an ideal of the algebra $A(S)^-.$ Since Witt algebra is simple, $J$ should be trivial. So, the algebra  $A(S)$ is simple. 

These arguments show  that classes of Leibniz algebras and weak Leibniz algebras are different. 
We do not know whether exists simple weak Leibniz algebras except algebras of a form $A(S).$

\subsection{Dialgebas, polarizations and depolarizations} 
A dialgebra is a linear space with two binary operations.  
A dialgebra $(A,\circ,\bullet)$ is called {\it polarization} of algebra $(A,\times)$ if 
$$a\circ b=\frac{1}{2}(a\times b-b\times a),$$
$$a\bullet b=\frac{1}{2}(a\times b+b\times a),$$
for any $a,b\in A.$ Conversely, an algebra $(A,\times)$ is called {\it depolarization} of 
dialgebra  $(A,\circ,\bullet)$ if
$$a\times b=a\circ b+a\bullet b$$
for any  $a,b\in A.$

A subspace $J\subseteq A$ is called {\it ideal}  of dialgebra $(A,\circ,\bullet),$ if 
$J$ is ideal of algebra $(A,\circ),$ 
$$a\circ b\in J,\quad  b\circ a\in J,$$ 
for any $a\in J, b\in A$ and  $J$ is ideal of algebra $(A,\bullet),$ 
$$a\bullet b\in J, \quad b\bullet a\in J,$$ 
for any $a\in J, b\in A.$ It is clear that ideal of dialgebra  $(A,\circ,\bullet)$ generates an ideal of its depolarization $(A,\times ).$ Conversely, an ideal of algebra $(A,\times)$ generates an ideal of its polarization $(A,\circ,\bullet).$ Therefore, the following statement holds.

\begin{prp} A dialgebra $(A,\circ, \bullet)$ is simple if and only if its polarization $(A,\times)$ is simple. 
\end{prp}

\subsection{Transposed Poisson algebras} 
A dialgebra $(A,\circ ,\bullet)$ is called {\it transposed Poisson } \cite{Bai}  if 
\begin{itemize}
\item an algebra $(A,\circ)$ is Lie
\item an algebra $(A,\bullet)$ is associative commutative, and 
\item associative part acts on Lie part as $1/2$-derivation,
$$2\,a\bullet(b\circ c)-(a\bullet b)\circ c-b\circ (a\bullet c)=0,$$
for any $a,b,c\in A.$
\end{itemize}

\subsection{Main results} 
\begin{thm} \label{mainbir} If $(A,\times)$ is weak Leibniz over a field of characteristic $p\ne 2,3,$ then its polarization   $(A,\circ,\bullet)$ is transposed Poisson. 

Conversely, if $(A,\circ,\bullet)$ is transposed Poisson algebra over a field of characteristic $p\ne 2,$ then 
its depolarization $(A,\times)$ is weak Leibniz. 
\end{thm}

\begin{crl} $(p\ne 2,3)$ An algebra is weak Leibniz iff its polarization is transposed Poisson. 
\end{crl}

\begin{thm} \label{dual} $(p\ne 2,3)$ Weak Leibniz operad is self-dual, i.e., its Koszul dual is equivalent to weak Leibniz operad.
Weak Leibniz operad is not Koszul. 
\end{thm}

{\bf Remark.} Any weak Leibniz algebra is associative-admissible and two-sided Alia, but not Poisson-admissible. In particular, any weak Leibniz algebra is Lie-admissible and 1-Alia.  On Alia algebras see \cite{Dzh2}. 

\section{An algebra is weak Leibniz iff its polarization is transposed Poisson} 

Let $K[t_1,t_2,\ldots ]$ be free magmatic algebra with generators $t_1,t_2,\ldots .$ 
Sometimes free magmatic algebra are called absolute free or non-commutative non-associative. This means that for any fixed $n$ there are $c_n=\frac{1}{n+1}{2n\choose n}$ kinds of bracketing types and any non-commutative non-associative polynomal is a combination of $c_n$ monomials. One can imagine these monomials as binary planar rooted trees with $n$ leafes, where any inner vertex is interpretered as a product of elements corresponding to  its sons. For example, monomials constructed by four generators $a,b,c,d$ are
$$
\begin{picture}(200,80)
\put(-50,40){$c_4=5$}
\put(30,30){\line(-1,1){30}}
\put(30,30){\line(1,1){10}}
\put(10,50){\line(1,1){10}}
\put(20,40){\line(1,1){10}}
\put(0,10){$((ab)c)d$}

\put(80,30){\line(-1,1){20}}
\put(80,30){\line(1,1){10}}
\put(80,50){\line(-1,1){10}}
\put(70,40){\line(1,1){20}}
\put(55,10){$(a(bc))d$}

\put(130,30){\line(-1,1){25}}
\put(130,30){\line(1,1){25}}
\put(115,45){\line(1,1){10}}
\put(145,45){\line(-1,1){10}}
\put(110,10){$(ab)(cd)$}

\put(180,30){\line(-1,1){10}}
\put(180,30){\line(1,1){20}}
\put(190,40){\line(-1,1){20}}
\put(180,50){\line(1,1){10}}
\put(165,10){$a((bc)d)$}

\put(230,30){\line(-1,1){10}}
\put(230,30){\line(1,1){30}}
\put(240,40){\line(-1,1){10}}
\put(250,50){\line(-1,1){10}}
\put(220,10){$a(b(cd))$}

\sudda{
\put(150,30){\line(-1,1){30}}
\put(150,30){\line(1,1){30}}
\put(130,50){\line(1,1){10}}
\put(170,50){\line(-1,1){10}}
}
\end{picture}
$$

For an algebra $A$ and for a polynomial $f=f(t_1,\ldots,t_k)\in K[t_1,t_2,\ldots ]$ we can make substitutions generators $t_i$ by elements of algebra $a_i\in A$ and calculate $f(a_1,a_2,\ldots,a_k)\in A$ in terms of multipliactions of $A.$ We say that $A$ satisfies  {\it polynomial identity } $f=0$ if 
$$f(a_1,a_2,\ldots,a_k)=0,\quad \forall a_1,a_2,\ldots,a_k\in A.$$
A class of algebras that satisfy identities $f_1=0,\ldots,f_s=0$ is called {\it variety } of algebras and is denoted $\langle f_1,\ldots,f_s\rangle.$

For polynomials $f_1,\ldots,f_s, g\in K[t_1,t_2,\ldots]$ we say that $g=0$ is {\it a consequence of}  identities $f_1=0,\ldots,f_s=0$  and denote $\{f_1, \ldots,f_s\}\Rightarrow g,$ if any algebra $A\in\langle f_1,\ldots,f_s\rangle$ satisfies the identity $g=0$ also. Easy generalization of this notion $\{f_1,\ldots,f_s\}\Rightarrow \{g_1,\ldots,g_l\}$ is clear. This means that any algebra
$A\in\langle f_1,\ldots,f_s\rangle$ satisfies identities $g_1=0,\ldots g_l=0$ also.  
In other words, notation $\{f_1,\ldots,f_s\}\Rightarrow \{g_1,\ldots,g_l\}$ means that
$\langle f_1,\ldots,f_s\rangle \subseteq 
\langle g_1,\ldots,g_l\rangle.$

We say that polynomial identities $f_1=0,\ldots, f_s=0$ and $g_1=0,\ldots,g_l=0$ are equivalent and write $\{f_1,\ldots,f_s\}\sim \{g_1,\ldots,g_l\},$ if 
$$\{f_1,\ldots,f_s\}\Rightarrow \{g_1,\ldots,g_l\} \mbox{ and }
\{g_1,\ldots,g_l\}\Rightarrow \{f_1,\ldots,f_s\}.$$

 Let us define polynomials $lwlei$ (left-weak Leibniz), $rwlei$ (right-weak Leibniz), $lieadm$ (Lie-admissibile), $assadm$ (associative-admissible)   and $1/2$-derivation polynomial or "associative part acts on Lie part" polynomial  $alder$ by 
 $$lwlei=[t_1,t_2]t_3-2t_1(t_2t_3)+2t_2(t_1t_3),$$
 $$rwlei=t_1[t_2,t_3]-2(t_1t_2)t_3+2(t_1t_3) t_2,$$ 
 $$lieadm=[[t_1,t_2],t_3]+[[t_2,t_3],t_1]+[[t_3,t_1],t_2],$$
 $$assadm=\{t_1,\{t_2,t_3\}\}-\{\{t_1,t_2\},t_3\},$$
$$alder(t_1,t_2,t_3)=2\{[t_1,t_2],t_3\}-[t_1,\{t_2,t_3\}]-[\{t_1,t_3\},t_2],$$
where
$$[t_1,t_2]=t_1t_2-t_2t_1,$$
$$\{t_1,t_2\}=t_1t_2+t_2t_1.$$

Let us establish relations between  polynomial identities generated by these polynomials.

\begin{lm}\label{40}
 $$assadm(t_1,t_2,t_3)+lieadm(t_1,t_2,t_3)-lwlei(t_3,t_1,t_2)+rwlei(t_2,t_3,t_1)=0.$$
\end{lm}

{\bf Proof.}  We have 
 $$assadm(t_1,t_2,t_3)+lieadm(t_1,t_2,t_3)=$$
 $$\{t_1,\{t_2,t_3\}\}-\{\{t_1,t_2\},t_3\}+
 [[t_1,t_2],t_3]+ [[t_2,t_3],t_1]+ [[t_3,t_1],t_2]=$$
 
$$2 t_1 (t_3 t_2) - 2 t_3 (t_1 t_2) - 2 (t_2 t_1) t_3 + 
 2 (t_2 t_3) t_1+$$
$$t_2 (t_1 t_3) - t_2 (t_3 t_1) - (t_1 t_3) t_2+ 
 (t_3 t_1) t_2=$$
 
$$- (t_1 t_3) t_2+ (t_3 t_1) t_2
 +2 t_1 (t_3 t_2) - 2 t_3 (t_1 t_2) -$$
  $$ 2 (t_2 t_1) t_3 + 
 2 (t_2 t_3) t_1+t_2 (t_1 t_3) - t_2 (t_3 t_1)  =$$
 
 $$lwlei(t_3,t_1,t_2)-rwlei(t_2,t_3,t_1).$$
 $\square$

\begin{lm}\label{4} If $p\ne 2,$ then 
$$
 \{lwlei, rwlei\} \sim \{lwlei, alder\}\sim \{ rwlei, alder\}
$$
 \end{lm}

{\bf Proof.}
We have 
$$alder(t_1,t_2,t_3)=lwlei(t_1,t_3,t_2)-lwlei(t_2,t_3,t_1)-rwlei(t_1,t_2,t_3)+rwlei(t_2,t_1,t_3).$$
Therefore, 
$$ \{lwlei, rwlei\} \Rightarrow \{lwlei, alder\}$$

Further, 
$$rwlei(t_1,t_2,t_3)=-1/2\,alder(t_1,t_2,t_3)+1/2\,alder(t_1,t_3,t_2)+1/2\,alder(t_2,t_3,t_1)-lwlei(t_2,t_3,t_1).$$
Therefore,
$$ \{lwlei, alder\} \Rightarrow \{lwlei, rwlei\}.$$

Hence,
$$ \{lwlei, rwlei\} \sim \{lwlei, alder\}.$$

The equivalence $\{lwlei, rwlei\} \sim \{ rwlei, alder\}$ is established by similar calculations.
$\square$

\begin{lm} \label{51} If $p\ne 3,$ then $ \{lwlei, rwlei\} \Rightarrow  \{lieadm, assadm\}$
\end{lm}

{\bf Proof.} Let $lalia$ and $ralia$ are left-alia and right-alia polynomials \cite{Dzh2},
$$lalia(t_1,t_2,t_3)=[t_1,t_2]t_3+[t_2,t_3]t_1+[t_3,t_1]t_2,$$
$$ralia(t_1,t_2,t_3)=t_1[t_2,t_3]+t_2[t_3,t_1]+t_3[t_1,t_2].$$

Note that the polynomial $lwlei(t_1,t_2,t_3)$ is skew-symmetric by 
variables $t_1,t_2,$ and the polynomial $rwlei(t_1,t_2,t_3)$ is skew-symmetric by variables $t_2,t_3,$ 
$$lwlei(t_1,t_2,t_3)+lwlei(t_2,t_1,t_3)=0,$$
$$rwlei(t_1,t_2,t_3)+rwlei(t_1,t_3,t_2)=0.$$
Therefore, the polynomial 
$$L(t_1,t_2,t_3)=lwlei(t_1, t_2, t_3) +lwlei(t_3, t_1, t_2) + 
   lwlei(t_2, t_3, t_1)$$   
is skew-symmetric by all variables $t_1,t_2,t_3$ and $L(t_1,t_2,t_3)$   is a skew-symmetrization of the term $(t_1t_2)t_3-2t_1(t_2t_3).$ We see that
$$L=lalia-2\,ralia,$$
By similar reasons, the polynomial
$$R(t_1,t_2,t_3)=
rwlei(t_1, t_2, t_3) +rwlei(t_3, t_1, t_2) + 
   rwlei(t_2, t_3, t_1)$$
is skew-symmetric by all variables $t_1,t_2,t_3,$ and
$$R=ralia-2\,lalia.$$
Therefore
$$L-R=3(lalia-ralia).$$
Since 
$$lieadm=lalia-ralia,$$
for  $p\ne 3$ we  have
$$lieadm(t_1,t_2,t_3)=$$ $$
1/3 \{lwlei(t_1, t_2, t_3) + lwlei(t_3, t_1, t_2) + 
   lwlei(t_2, t_3, t_1) - rwlei(t_1, t_2, t_3) - $$ $$
   rwlei(t_2, t_3, t_1) - rwlei(t_3, t_1, t_2)\}.$$
 In other words,
 $$\{lwlei,rwlei\}\Rightarrow lieadm.$$
 Hence in case $p\ne 3$ by Lemma \ref{40},
  $$\{lwlei,rwlei\}\Rightarrow assadm.$$

  $\square$

 \begin{lm}\label{52} $ \{lwlei, rwlei\} \Rightarrow  alder$
\end{lm}
 
 {\bf Proof.} 
%For$$alder(t_1,t_2,t_3)=2\{[t_1,t_2],t_3\}-[t_1,\{t_2,t_3\}]+[t_2,\{t_1,t_3\}]$$
We have 
$$alder(t_1,t_2,t_3)=lwlei(t_1, t_3, t_2) - lwlei(t_2, t_3, t_1) - 
 rwlei(t_1, t_2, t_3) + rwlei(t_2, t_1, t_3).$$
So,
$$ \{lwlei, rwlei\} \Rightarrow  alder.$$
$\square$

\begin{lm} \label{53} If $p\ne 2,$, then $\{lieadm,assadm, alder\}\Rightarrow \{lwlei,rwlei\}$
\end{lm}

{\bf Proof.}  We have 
$$lwlei(t_1,t_2,t_3)=$$
  $$ 1/4 \{2 \,assadm(t_2, t_3, t_1) + 2\, lieadm(t_1, t_2, t_3) + 
   alder(t_1, t_2, t_3) + aldder(t_1, t_3, t_2)$$ $$ - 
   alder(t_2, t_3, t_1)\},$$
   
\sudda{$$rwlei(t_1,t_2,t_3)=$$
 $$  1/4 \{2\, assadm(t_1, t_2, t_3) + 2\; assadm(t_2, t_3, t_1) - 
   2 \,lieadm(t_1, t_2, t_3) - alder(t_1, t_2, t_3) + $$ $$
   alder(t_1, t_3, t_2) + alder(t_2, t_3, t_1)\}.$$
   }
   
So, $$ \{lieadm, assadm, alder\}\Rightarrow \{lwlei, alder\} , $$
and by Lemma \ref{4}
$$ \{lieadm, assadm, alder\}\Rightarrow \{lwlei, rwlei\} $$
$\square$

{\bf Proof of Theorem \ref{mainbir}.} If $p\ne 3,$ by Lemmas \ref{51} and \ref{52}
$$ \{lwlei, rwlei\} \Rightarrow  \{lieadm, assadm,alder\} .$$
If $p\ne 2,$ by Lemma \ref{53}
$$ \{lieadm, assadm,alder\} \Rightarrow \{lwlei, rwlei\} .$$
So, if $p\ne 2 ,3,$ 
 $$\{lieadm, assadm,alder\} \sim \{lwlei, rwlei\} $$
 $\square$

\section{Novikov-Poisson algebras and transposed Poisson algebras}

\subsection{Novikov-Poisson algebras} 
An algebra $(A,\cdot,\bullet)$ is called {\it Novikov-Poisson,} \cite{Xu} if 
\begin{itemize}
\item $(A,\cdot)$  is (left) Novikov, for any $a,b,c\in A,$
$$(a\cdot b -b\cdot a)\cdot c=a\cdot (b\cdot c)-b\cdot (a\cdot c), \quad (a\cdot b)\cdot c=(a\cdot c)\cdot b,$$
\item $(A,\bullet)$ is associative commutative, such that for any $a,b,c\in A,$
$$a\bullet (b\cdot c)=(a\bullet b)\cdot c, \quad a\cdot (b\bullet c)=(a\cdot b)\bullet c+b\bullet (a\cdot c).$$
\end{itemize}

\subsection{Construction of transposed Poisson algebras by Novikov-Poisson algebras}

\begin{prp}\label{A_{u,v}}
 Let $A=(A,\cdot,\bullet)$ be Novikov-Poisson algebra.  Then for any $u,v\in A$ the algebra $A_{u,v}=(A,\circ_u,\bullet_v),$ where
$$a\circ_u v=u\bullet(a\cdot b-b\cdot a), \quad a\bullet_v b=v\bullet(a\bullet b),$$
is transposed Posson and the algebra $A_{u,v}$ under multiplication $a b=1/2(a\circ_u b+a\bullet_v b)$ is weak Leibniz. 
\end{prp}

Proof of this proposition is based on  the following two lemmas.

\begin{lm} \label{2022bir}
Let $A=(A,\cdot,\bullet)$ is Novikov-Poisson, $u\in A$ and 
$$a\cdot_u b=u\bullet(a\cdot b).$$
Then $A_u=(A,\cdot_u,\bullet ) $ is Novikov-Poisson for any $u\in A.$
\end{lm}

{\bf Proof.} Since $A=(A,\cdot,\bullet)$ is Novikov-Poisson, by definition 
the algebra $(A,\bullet)$ is associative-commutative. 

 Let us check that $(A,\cdot_u)$ is left-Novikov,
Let $$(a,b,c)_u=a\cdot_u(b\cdot_u c)-(a\cdot_u b)\cdot_u c.$$
Then 
$$(a,b,c)_u=$$

$$u\bullet(a\cdot (u\bullet(b\cdot c)))-u\bullet((u\bullet(a\cdot b))\cdot c)=$$

$$u\bullet(a\cdot (u)\bullet(b\cdot c)))+
u\bullet(u\bullet(a\cdot (b\cdot c)))
-(u\bullet u)\bullet((a\cdot b)\cdot c)=$$

$$(a\cdot (u)\bullet u)\bullet(b\cdot c)+
(u\bullet u)\bullet(a\cdot (b\cdot c)-(a\cdot b)\cdot c)=$$

$$1/2(a\cdot (u\bullet u))\bullet(b\cdot c)+(u\bullet u)\bullet (a,b,c).$$

Similarly,
$$(b,a,c)_u=$$
$$1/2(b\cdot (u\bullet u))\bullet(a\cdot c)+(u\bullet u)\bullet (b,a,c).$$
Hence
$$lsym_u(a,b,c)=(a,b,c)_u-(b,a,c)_u=$$
$$1/2(a\cdot (u\bullet u))\bullet(b\cdot c)-
1/2(b\cdot (u\bullet u))\bullet(a\cdot c)=$$

$$1/2\{
(b\cdot c)\bullet (a\cdot (u\bullet u))-
(a\cdot c)\bullet (b\cdot (u\bullet u))\}=$$

$$1/2\{(b\cdot c)\bullet a-
(a\cdot c)\bullet b\}\cdot (u\bullet u)=$$

$$1/2\{a\bullet (b\cdot c)-
b\bullet (a\cdot c)\}\cdot (u\bullet u)=$$

$$1/2\{(a\bullet b)\cdot c)-
(b\bullet a)\cdot c\}\cdot (u\bullet u)=$$

$$0.$$

We have
$$a\cdot_u (b\bullet c)=u\bullet (a\cdot(b\bullet c))=$$
$$u\bullet ((a\cdot b)\bullet c+b\bullet (a\cdot c))=$$
$$(u\bullet (a\cdot b))\bullet c+u\bullet (b\bullet (a\cdot c))=$$
$$(u\bullet (a\cdot b))\bullet c+ b\bullet (u\bullet(a\cdot c))=$$
$$(a\cdot_u b)\bullet c+b\bullet (a\cdot_u c).$$
Further,
$$(a\bullet b)\cdot_u c=u\bullet ((a\bullet b)\cdot c)=
u\bullet (a\bullet (b\cdot c))=
 a\bullet (u\bullet(b\cdot c))=
a\bullet (b\cdot_u c).$$
$\square$

\begin{lm} \label{2022eki} Let $A=(A,\cdot,\bullet)$ is Novikov-Poisson, $u\in A,$  and 
$$a\bullet_u b=u\bullet(a\bullet  b).$$
Then $A_u=(A,\cdot,\bullet_u ) $ is transposed Poisson for any $u\in A.$
\end{lm}

{\bf Proof.}  Since $A=(A,\cdot,\bullet)$ is Novikov-Poisson, the algebra
$(A\cdot)$ is Lie and the algebra $(A,\bullet_u)$ is associative commutative. 

Let us check $1/2$-derivation condition. For $a\circ b=a\cdot b-b\cdot a$ we have 
$$(a\bullet_u b)\circ c=(u\bullet (a\bullet b))\circ c =(a\bullet b \bullet u)\cdot c-c\cdot(a\bullet b\bullet  u)=$$
$$(a\bullet b \bullet_u)\cdot c-a\bullet b\bullet  c\cdot(u)
-a\bullet c\cdot(b)\bullet u-c\cdot(a)\bullet b\bullet  u,$$

$$b\circ (a\bullet_u c)=b\cdot (a\bullet_u c)- (a\bullet_u c)\cdot b=
b\cdot (u\bullet a\bullet c)-( u\bullet a\bullet c)\cdot b=$$

$$(b\cdot u)\bullet a\bullet c)+
u\bullet (b\cdot a)\bullet c)+
u\bullet a\bullet (b\cdot c)
- u\bullet a\bullet (c\cdot b).$$
Hence,
$$(a\bullet_u b)\circ c+b\circ (a\bullet_u c)=$$

$$(a\bullet b \bullet u)\cdot c-a\bullet b\bullet ( c\cdot u)
-a\bullet (c\cdot b)\bullet u-(c\cdot a)\bullet b\bullet  u+$$
$$(b\cdot u)\bullet a\bullet c+
u\bullet (b\cdot a)\bullet c+
u\bullet a\bullet (b\cdot c)
- u\bullet a\bullet (c\cdot b)=$$

$$(a\bullet b \bullet u)\cdot c+u\bullet a\bullet (b\cdot c)+$$
$$-a\bullet b\bullet ( c\cdot u)+(b\cdot u)\bullet a\bullet c$$
$$-a\bullet (c\cdot b)\bullet u- u\bullet a\bullet (c\cdot b)$$
$$-(c\cdot a)\bullet b\bullet  u+u\bullet (b\cdot a)\bullet c=$$

$$(a\bullet b \bullet u)\cdot c+(u\bullet a\bullet b)\cdot c+$$
$$-(a\bullet b\bullet  c)\cdot u)+ (a\bullet c\bullet b)\cdot u$$
$$-(a\bullet u\bullet c)\cdot b- (u\bullet a\bullet c)\cdot b)$$
$$-(b\bullet  u\bullet c)\cdot a+(u\bullet c \bullet b)\cdot a=$$

$$2(a\bullet b \bullet u)\cdot c- 2(u\bullet a\bullet c)\cdot b)=$$
Further,
$$2a\bullet_u (b\circ c)=u\bullet (2a\bullet(b\circ c))=2 (a\bullet b\bullet u)\cdot c-2 (a\bullet c\bullet u)\cdot b.$$
Therefore,
$$2a\bullet_u (b\circ c)=(a\bullet_u b)\circ c+b\circ (a\bullet_u c).$$
$\square$

{\bf Proof of Proposition \ref{A_{u,v}}.}
\sudda{
\begin{crl} \label{2022ush}
 Let $A=(A,\cdot,\bullet)$ is Novikov-Poisson, $u,v\in A,$  and 
$$a\cdot_u b=u\bullet (a\cdot b),$$
$$a\bullet_v b=v\bullet(a\bullet  b).$$
Then $A_{u,v}=(A,\cdot_u,\bullet_v ) $ is transposed Poisson for any $u,v\in A.$
\end{crl}
}
By Lemma \ref{2022bir} the algebra 
$A_u=(A,\cdot_u,\bullet)$ is Novikov-Poisson. Therefore, by Lemma \ref{2022eki}  the algebra $A_{u,v}=(A,\cdot_u,\bullet_v)$ is transposed Poisson. 
$\square$

\section{Weak Leibniz algebra is self-dual and is not Koszul}

{\bf Proof of Theorem \ref{dual}.} Base of weak Leibniz algebras
$$[a,b]c-2a(bc)+2b(ac)=0$$
$$a[b,c]-2(ab)c+2(ac)b=0$$
in degree 3 has dimension $6$ and base elements are 
$$\{ c(a b), (bc)a, (ca)b, (ac)b, (ba)c,(ab)c \}.$$
 
Other non-base elements are 

$$a (c b)= c (a b)+ 1/2 (a c) b- 1/2 (c a) b,$$
  
$$ a (b c)=  c (a b)+ 2 (a b) c- 3/2 (a c) b-    1/2 (c a) b,$$

$$b (c a)= c (a b) + 3/2 (a b) c - 3/2 (a c) b- 3/2 (b a) c+ 2 (b c) a - 1/2 (c a) b,$$
 
$$ b (a c)= c (a b) + 3/2 (a b) c - 3/2 (a c) b +    1/2 (b a) c- 1/2 (c a) b ,$$
   
$$c (b a)= c (a b) + 2 (a b) c- 2 (a c) b -  2(b a) c + 2 (b c) a,$$
    
$$ (c b) a= (a b) c - (a c) b - (b a) c + (b c) a + (c a) b.$$

To find Koszul dual we have to 
calculate Jacobi polynomial for $A\otimes U,$ where $A$ is weak Leibniz and $U$ is its Koszul dual. 
We have 
$$jac(a\otimes u,b\otimes v,c\otimes w)=$$

 $$=c (a b) \otimes (-u (v w)+u (w v)+v (u w)- v (w u)-w (u v)+w (v u))$$
 
$$+(b c) a\otimes  (-2 (v (w u)+ 2 (w (v u)+(v w) u- (w v) u))$$
 
$$+ 1/2 (c a)b\otimes ((u (v w)-u (w v)-v (u w)+     v (w u)+ 2 (w u) v- 2 (w v) u)$$
 
$$+ 1/2 (a c)  b\otimes  (3 u (v w)+u (w v)- 3 v (u w)+     3 v (w u)- 4 w (v u)- 2 (u w) v+ 
    2 (w v) u))$$
  
$$+ 1/2 (b a) c\otimes (v (u w)+ 3 v (w u)- 4 w (v u)- 
    2 (v u) w+ 2 (w v) u))$$

$$+ 1/2 (a b) c \otimes (-4 u (v w)+ 3 v (u w)- 3 v (w u)+ 
    4w (v u)% here was 4, but should be 2 ??
    + 2 (u v) w- 2 (w v) u))
  $$

We see that  
$$-u (v w)+u (w v)+v (u w)- v (w u)-w (u v)+w (v u)=-u[v,w]-v[w,u]-w[u,v]=0,$$

$$ (-2 (v (w u)+ 2 (w (v u)+(v w) u- (w v) u))=0,$$

$$((u (v w)-u (w v)-v (u w)+     v (w u)+ 2 (w u) v- 2 (w v) u)=0,$$

$$(3 u (v w)+u (w v)- 3 v (u w)+     3 v (w u)- 4 w (v u)- 2 (u w) v+ 
    2 (w v) u))=0,$$
    
$$ (v (u w)+ 3 v (w u)- 4 w (v u)-
    2 (v u) w+ 2 (w v) u))=0,$$
    
$$ (-4 u (v w)+ 3 v (u w)- 3 v (w u)+ 
    4 w (v u)% 4 w (v u) but should be 2??
      + 2 (u v) w- 2 (w v) u))=0.$$

Conversely, these identities imply weak Leibniz identities.
Hence dual operad for weak Leibniz operad is generated by identities
$$u[v,w]-2(uv)w+2(uw)v=0,$$
$$[u,v]w-2 u (v w)+ 2 v(uw)=0.$$
that are weak Leibniz identities.

For weak Leibniz operads dimensions  of multi-linear parts are as below
$$\begin{array}{|c|c|c|c|c|c|c|}
\hline
n&1&2&3&4&5&6\\
\hline
d_n=d_n^!&1&2&6&20&74&301\\
\hline
\end{array}
$$
Hence beginning parts of skew-exponential generating functions are
$$f(x)=f^!(x)=-x + x^2 - x^3 + (5 x^4)/6 - (37 x^5)/60,$$
and
$$f(f^!(x))=x+7x^5/30+O(x)^6\ne x.$$
Thus by Koszulity necessary condition \cite{GK} the weak Leibniz operad is not Koszul.

\end{document}